\documentclass[12pt]{elsart}
\usepackage{amsfonts,amscd,amsmath,amssymb,pb-diagram,graphicx,latexsym}

\renewcommand{\L}{\Lambda}
\newcommand{\Aut}{\operatorname{Aut}}
\newcommand{\Inn}{\operatorname{Inn}}
\renewcommand{\O}{{\mathcal O}}
\newcommand{\Hom}{\operatorname{Hom}}
\newcommand{\im}{\operatorname{im}}
\newcommand{\vphi}{\varphi}
\newcommand{\diag}{\operatorname{diag}}
\newcommand{\Z}{\mathbb{Z}}
\newcommand{\rad}{\operatorname{rad}}

\begin{document}
\begin{frontmatter}
\title{Automorphisms of tiled orders}

\author{Jeremy Haefner}
\address{Department of Mathematics, University of Colorado, Colorado Springs, CO 80933-7150}

\author{Christopher J. Pappacena\thanksref{Baylor}}
\address{Department of Mathematics, Baylor University, Waco, TX 76798}

\thanks[Baylor]{Supported in part by a University Research Grant and a Summer Sabbatical from Baylor University.}

\begin{abstract} Let $\L$ be a tiled $R$-order.  We give a description of $\Aut_R(\L)$ as the semidirect product of $\Inn(\L)$ and a certain subgroup of $\Aut(Q(\L))$, where $Q(\L)$  is the link graph of $\L$.  Additionally, we give criteria for determining when an element of $\Aut(Q(\L))$ belongs to this subgroup in terms of the exponent matrix for $\L$.          
\end{abstract}


\end{frontmatter}

\section{Introduction}
Let $R$ be a complete discrete valuation ring (DVR) with quotient field $K$ and maximal ideal $P=\pi R$.  Recall that an $R$-order $\L$ in $M_n(K)$ is called \emph{tiled} in case it has a complete set of $n$ orthogonal idempotents $\{e_1,\dots, e_n\}$.  Without loss of generality, we may assume that $e_i=e_{ii}$ for all $i$, where $\{e_{ij}\}$ are the standard matrix units in $M_n(K)$.  It is then possible to write $\L=(P^{\alpha_{ij}})_n=(\pi^{\alpha_{ij}}R)_n$, where the $\alpha_{ij}$ are non-negative integers, $\alpha_{ii}=0$ for all $i$, and $\alpha_{ij}+\alpha_{jk}\geq\alpha_{ik}$ for all $i$, $j$ and $k$. \cite[p. 313]{Jat}.  If $\Lambda$ is such a tiled order, then the integral matrix $(\alpha_{ij})$ is called the \emph{exponent matrix} for $\L$. It is clear that $\L$ is completely determined by its exponent matrix.

In this paper, we study the structure of $\Aut_R(\L)$, the group of $R$-automorphisms of $\L$, in the case where $\L$ is tiled.  We give a description of $\Aut_R(\L)$ as the  semidirect product of $\Inn(\L)$ and a certain subgroup of the automorphsism group of the \emph{link graph} $Q(\L)$ of $\L$.  (The precise definition is recalled in the next section; for now we note only that $Q(\L)$ is a quiver on $n$ vertices.)  This subgroup, which we denote by $\O_\L$, is called the group of \emph{liftable} automorphisms of $Q(\L)$.  As the name suggests, these automorphisms are precisely those which are induced by automorphisms of the order $\L$.  

In general, $\O_\L$ will be a proper subgroup of $\Aut(Q(\L))$; we develop necessary and sufficient conditions for a given automorphism to be liftable in terms of the matrix of exponents $(\alpha_{ij})$ of $\L$.  Additionally, we show that if $\L$ is basic and $\alpha_{ij}\in\{0,1\}$ for all $i,j$, then $\O_\L=\Aut(Q(\L))$.  Finally, using these liftable automorphisms, we give an example to show that the crossed product $\L*\O_\L$ need not be a prime ring.  This is in contrast to the hereditary case \cite{Haef Jan}.

\section{The link graph of a tiled order}
Let $\L=(P^{\alpha_{ij}})_n$ be a tiled order as above.  In \cite[section 2]{Wie Rog},  Wiedemann and Roggenkamp construct a quiver $Q(\L)$ associated to $\L$, as follows. Let $P_i=\Lambda e_i$ denote the $i$-th column of $\Lambda$.  Then $Q(\L)$ has $n$ vertices, and there is an arrow from $i$ to $j$ if and only if $P_i$ is a summand of the projective cover of $\rad(P_j)$. One assigns a value $v$ to the arrows by setting $v(i,j)=\beta$, where $\Hom_\L(P_i,P_j)=P^\beta$.  The corresponding valued quiver will be denoted $Q^v(\L)$.  In general, distinct orders $\L$ and $\L'$ can have identical quivers: $Q(\L)=Q(\L')$ (see for example \cite[Example 10]{Fujita Link}).  However, it is proven in \cite[Theorem 1]{Wie Rog} that $\L$ is uniquely determined by $Q^v(\L)$.  

There is an alternative description of $Q(\L)$ as the \emph{link graph} of the maximal ideals of $\L$, which we now describe.  Note that $\L$ has exactly $n$ distinct maximal two-sided ideals $M_1,\dots, M_n$, where $M_k$ is obtained from $\L$ by repacing the $R$ in the $(k,k)$-position with a $P$. That is, $M_k=(P^{\beta_{ij}})_n$, where $\beta_{ij}=\alpha_{ij}$ for $(i,j)\neq (k,k)$, and $\beta_{kk}=1$. The \emph{link graph} is defined as the quiver with $n$ vertices, with an arrow from $i$ to $j$ if and only if $M_jM_i\neq M_j\cap M_i$ \cite{Muller FBN}. Then, \cite[Proposition 1.2]{Fujita rem} shows that the link graph of $\L$ coincides with $Q(\L)$. 

\begin{lem}  There is a group homomorphism $\Phi:\Aut_R(\Lambda)\rightarrow \Aut(Q(\L))$, whose kernel contains $\Inn(\L)$.\label{Phi lemma}
\end{lem}
\begin{pf}  An automorphism of $\L$ permutes the maximal ideals $M_1,\dots, M_n$ of $\L$, and thereby gives rise to a permutation of the vertices of $Q(\Lambda)$. Let us denote the corresponding permutation of $\{1,\dots, n\}$ as $\sigma$.  Then, we have that there is an arrow from $\sigma(i)$ to $\sigma(j)$ in $Q(\L)$ if and only if $M_{\sigma(j)}M_{\sigma(i)}\neq M_{\sigma(j)}\cap M_{\sigma(i)}$, if and only if $M_jM_i\neq M_j\cap M_i$, if and only if there is an arrow from $i$ to $j$ in $Q(\L)$.  Thus, $\sigma$ gives an automorphism of $Q(\L)$.  Finally, if $\vphi$ is an inner automorphism of $\L$, then $\vphi$ fixes each of the maximal ideals $M_1,\dots, M_n$.  Thus $\vphi$ induces the identity on $Q(\L)$.\qed
\end{pf}

We shall see below that in fact $\Inn(\L)=\ker\Phi$.  The reason that we work with link graph $Q(\L)$ instead of the valued quiver $Q^v(\L)$ is that the corresponding result fails for valued quivers: There are tiled orders $\L$ for which $\vphi\in\Aut(\L)$ does not induce an automorphism of $Q^v(\L)$. The following example illustrates this.

\begin{exmp}\normalfont Consider the order
\[\L=\begin{pmatrix} R&P^2&P^4\\P^3&R&P^4\\P&P&R\end{pmatrix}.\]
One computes that the valued quiver $Q^v(\L)$ is 

\centerline{\includegraphics{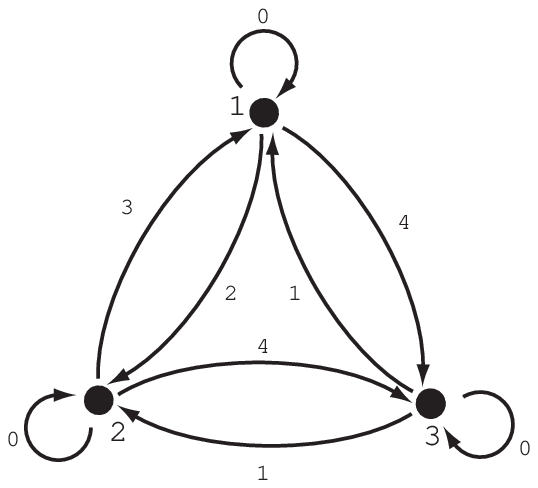}}

Now, a direct verification shows that conjugation by $\Bigl(\begin{smallmatrix} 0&\pi&0\\0&0&\pi^3\\1&0&0\end{smallmatrix}\Bigr)$ is an automorphism of $\L$, and the induced permutation on the vertices of $Q(\L)$ is $\sigma=(123)$.  Since $v(1,2)=2$ and $v(\sigma(1),\sigma(2))=v(2,3)=4$, we see that $\sigma$ does \emph{not} induce an automorphism of $Q^v(\L)$.\qed
\end{exmp}

One may at first hope that every automorphism of $Q(\L)$ is induced by an automorphism of $\Aut_R(\L)$ as above; i.e. that $\Phi$ is surjective.  Unfortunately, this fails to be the case.  Indeed, we shall develop an explicit criterion for determining whether or nor an element of $Q(\L)$ is in $\im \Phi$; we shall call such automorphisms \emph{liftable}. 

We first fix some notation.  Given $x\in GL_n(K)$, we denote conjugation by $x$ as $\iota_x$.  We shall identify automorphisms of $Q(\L)$ with the induced permutation of the vertices of $Q(\L)$; thus $\sigma$ will simultaneously denote an element of $\Aut(Q(\L))$ and an element of $S_n$.  Given $\sigma\in S_n$, we let $v(\sigma)\in M_n(K)$ denote the corresponding permutation matrix: $v(\sigma)_{ij}=1$ if $j=\sigma(i)$, $0$ otherwise.

\begin{lem} Let $\vphi\in\Aut_R(\L)$. Then there exist $u\in\L^*$, $d=\diag(\pi^{d_1},\dots,\pi^{d_n})$, and $\sigma\in S_n$ such that $\vphi=\iota_{udv(\sigma)}$.\label{factor automorphism}
\end{lem}

\begin{pf} The Skolem-Noether Theorem shows that any automorphism of $\L$ is given by conjugation by $x$ for some $x\in GL_n(K)$.  The automorphism $\vphi$ takes the set of orthogonal primitive idempotents $\{e_1,\dots,e_n\}$ to another set of orthogonal primitive idempotents, say $\{f_1,\dots,f_n\}$.  Since $\L$ is semiperfect, there is a $u\in\L^*$ and $\sigma\in S_n$ such that $\iota_uf_i=e_{\sigma(i)}$ for all $i$ \cite[Proposition 3.7.3]{Lambek}.  

So, the automorphism $\iota_u\vphi$ which is conjugation by some $x\in GL_n(K)$, acts as a permutation $\sigma$ on $\{e_1,\dots,e_n\}$.  If one writes out in explicit matrix form the condition that $x^{-1}e_ix=e_{\sigma(i)}$ for all $i$, then one sees that only the $(i,\sigma(i))$ entry of $x$ is nonzero, for $i=1,\dots, n$.  We may write the $(i,\sigma(i))$ entry of $x$ as $r_i\pi^{d_i}$, where $r_i\in R^*$ and $d_i\in \Z$.  In particular, $x$ factors as $ydv(\sigma)$, where $y=\diag(r_1,\dots,r_n)$ and  $d=\diag(\pi^{d_1},\dots,\pi^{d_n})$.  Now, since $\iota_u\vphi=\iota_{ydv(\sigma)}$, we have $\vphi=\iota_{u^{-1}ydv(\sigma)}$. Since each of $u^{-1}$ and $y$ is in $\L^*$, we see that $\vphi$ has the indicated form.\qed
\end{pf}

\begin{prop} An automorphism $\sigma$ of $Q(\L)$ is liftable if and only if there exists a diagonal matrix $d=\diag(\pi^{d_1},\dots,\pi^{d_n})$ such that $\iota_{dv(\sigma)}\in\Aut_R(\Lambda)$.\label{liftable lemma}
\end{prop}

\begin{pf}  Suppose that there is an automorphism of $\L$ of the given form.  Since conjugation by $d$ fixes the primitive orthogonal idempotents $\{e_1,\dots, e_n\}$, we see that $\iota_{dv(\sigma)}(e_i)=e_{\sigma(i)}$ for all $i$.  Since the maximal ideal $M_i$ of $\L$ can be characterized as the unique maximal ideal of $\L$ which does not contain $e_i$, this implies that $\iota_{dv(\sigma)}(M_i)=M_{\sigma(i)}$ for all $i$.  Thus $\iota_{dv(\sigma)}$ induces the permutation $\sigma$ on the vertices of $Q(\L)$, and is necessarily an automorphism of $Q(\L)$.  Hence, $\sigma=\Phi(\iota_{dv(\sigma)})$ is liftable.

Conversely, let $\sigma$ be liftable.  Then there is an automorphism $\vphi$ with $\Phi(\vphi)=\sigma$.  By Lemma \ref{factor automorphism}, we can write $\vphi=\iota_{udv(\tau)}$, where $u\in\Lambda^*$, $d=\diag(\pi^{d_1},\dots,\pi^{d_n})$, and $\tau\in S_n$.  Now, $\Phi(\vphi)=\Phi(\iota_{dv(\tau)})$ (since $\Phi(\iota_u)=0$) and, as in the previous paragraph, $\iota_{dv(\tau)}(M_i)=M_{\tau(i)}$ for all $i$.  It follows that $\tau=\sigma$, so that $\iota_{dv(\sigma)}\in\Aut_R(\L)$. \qed
\end{pf}

\section{Main results} We are now in a position to state and prove the main theorem of this paper.  

\begin{thm} Let $\L$ be a tiled order in $M_n(K)$ with exponent matrix $(\alpha_{ij})$ and link graph $Q(\L)$.  

\begin{enumerate}
\item[(a)] $\sigma\in\Aut(Q(\L))$ is liftable if and only if the linear system 
\begin{equation}
x_i-x_j=\alpha_{ij}-\alpha_{\sigma(i)\sigma(j)},\ \ i<j\label{system}\end{equation}
has a solution ${\bf x}=(x_1,\dots,x_n)$ in integers $x_i$.  

\item[(b)] If $\sigma$ is liftable, then $\iota_{d({\bf x})v(\sigma)}\in\Aut_R(\L)$, where $d({\bf x})=\diag(\pi^{x_1},\dots,\pi^{x_n})$.

\item[(c)] If $\vphi\in\Aut_R(\L)$ is written as $\iota_{udv(\sigma)}$ with  $d=\diag(\pi^{d_1},\dots,\pi^{d_n})$ as in Lemma \ref{factor automorphism}, then ${\bf x}=(d_1,\dots, d_n)$ is a solution to \eqref{system}. 

\item[(d)] Let $\O_\L=\{\iota_{d({\bf x})v(\sigma)}:\mbox{$\sigma$ is liftable}\}$. Then $\O_\L$ is a subgroup of $\Aut_R(\L)$, and $\Aut_R(\L)=\Inn(\L)\rtimes \O_\L$.
\end{enumerate}\label{main thm}
\end{thm}

\begin{pf} (a) Suppose that $\sigma$ is liftable. By Proposition \ref{liftable lemma}, $\sigma$ must lift to an automorphism of the form $\iota_{dv(\sigma)}$ for some diagonal matrix $d=\diag(\pi^{d_1},\dots, \pi^{d_n})$.  Now, the exponent matrix for the order $d^{-1}\L d$ is $(\alpha_{ij}-d_i+d_j)$; this follows because left multiplication by $d^{-1}$ subtracts $d_i$ from the $i$-th row of the exponent matrix, while right multiplication by $d$ adds $d_j$ to the $j$-th column.  Conjugating $d^{-1}\L d$ by $v(\sigma)$ has the effect of applying $\sigma$ to the indices in the exponent matrix; that is, the exponent matrix for $v(\sigma)^{-1}d^{-1}\L dv(\sigma)$ is $(\beta_{ij})=(\alpha_{\sigma(i)\sigma(j)}-d_{\sigma(i)}+d_{\sigma(j)})$.  Since $\iota_{dv(\sigma)}$ is an automorphism of $\L$, $\beta_{ij}=\alpha_{ij}$ for all $i$ and $j$; i.e.
\[
\alpha_{ij}=\alpha_{\sigma(i)\sigma(j)}-d_{\sigma(i)}+d_{\sigma(j)}
\]
for all $i$ and $j$.  It follows that ${\bf x}=(d_1,\dots, d_n)$ is a solution to \eqref{system}.

Conversely, let ${\bf x}=(x_1,\dots, x_n)$ be a solution to \eqref{system}, and let $d({\bf x})=\diag(\pi^{x_1},\dots, \pi^{x_n})$.  Then one computes as in the previous paragraph that the exponent matrix for $v(\sigma)^{-1}d({\bf x})^{-1}\L d({\bf x})v(\sigma)$ is $\alpha_{ij}$.  This shows that $\sigma$ is liftable.

(b)  The computation in the previous paragraph shows that $\iota_{d({\bf x})v(\sigma)}$ is an automorphism of $\L$ whenever $\sigma$ is liftable. 

(c)  Suppose $\vphi=\iota_{udv(\sigma)}$. Since $u\in\L^*$, we have that $\iota_{u^{-1}}\vphi=\iota_{dv(\sigma)}\in\Aut_R(\L)$.  Now, the computation in the proof of part (a) shows that ${\bf x}=(d_1,\dots, d_n)$ is a solution to \eqref{system}, where $d=\diag(\pi^{d_1},\dots, \pi^{d_n})$.  

(d) We first show that $\O_\L$ is a subgroup of $\Aut_R(\L)$.  Let $\iota_{dv(\sigma)}$ and $\iota_{\delta v(\tau)}$ be in $\Aut_R(\L)$, where $d=\diag(\pi^{d_1},\dots, \pi^{d_n})$ and $\delta=\diag(\pi^{\delta_1},\dots, \pi^{\delta_n})$.  Then $\iota_{dv(\sigma)}\iota_{\delta v(\tau)}=\iota_{dv(\sigma)\delta v(\tau)}$.  Now $dv(\sigma)\delta v(\tau)=d\delta^\sigma v(\sigma)v(\tau)=d\delta^\sigma v(\sigma\tau)$, where $\delta^\sigma=\diag(\pi^{\delta_{\sigma(1)}},\dots, \pi^{\delta_{\sigma(n)}})$.  Since $\iota_{d\delta^\sigma v(\sigma\tau)}$ induces $\sigma\tau$ on $Q(\L)$, we see that $d\delta^\sigma=d({\bf x})$ for some solution ${\bf x}$ of \eqref{system}, by part (c).  It follows that $\iota_{dv(\sigma)}\iota_{\delta v(\tau)}\in\O_\L$ and $\O_\L$ is a subgroup.   

Next, we show that $\O_\L\cap \Inn(\L)=1$.  If $\iota_{dv(\sigma)}\in\O_\L\cap\Inn(\L)$, then $\Phi(\iota_{dv(\sigma)})=1$ by Lemma \ref{Phi lemma}.  Thus $\sigma$ must be the identity automorphism of $Q(\L)$ by the proof of Proposition \ref{liftable lemma}.  Now, the system \eqref{system} has only the trivial solutions ${\bf x}=(x,\dots, x)$ for some $x\in \Z$, so that $d=\diag(\pi^x,\dots, \pi^x)$ for some $x\in\Z$. Thus $\iota_{dv(\sigma)}=\iota_d=1$.  

To finish the proof, we need to show that $\Aut_R(\L)=\Inn(\L)\O_\L$.  Given any $\vphi\in\Aut_R(\L)$, we can write $\vphi=\iota_u\iota_{dv(\sigma)}$ with $\iota_u\in\Inn(\L)$ by Lemma \ref{factor automorphism}.  By part (c) $\iota_{dv(\sigma)}\in\O_\L$, completing the proof. \qed
\end{pf}

\begin{exmp} \normalfont We illustrate Theorem \ref{main thm} for the order
\[\L=\begin{pmatrix} R&P^2&P^4\\P^3&R&P^4\\P&P&R\end{pmatrix}.\]
One computes that the link graph of $\L$ is the complete quiver on $3$ vertices, so that $\Aut(Q(\L))=S_3$.  If we let $\tau=(12)$ and $\sigma=(123)$, then $\tau$ and $\sigma$ generate $S_3$.

We first show that $\sigma$ is liftable.  The linear system of equations becomes
\begin{equation}
\begin{split}
x_1-x_2&=\alpha_{12}-\alpha_{23}=-2\\
x_1-x_3&=\alpha_{13}-\alpha_{21}=1\\
x_2-x_3&=\alpha_{23}-\alpha_{31}=3,\end{split}
\end{equation}
which has the solution $x_1=1+a$, $x_2=3+a$, $x_3=a$.  Taking $a=0$, we see that conjugation by $dv(\sigma)$ is an automorphism of $\L$, where $d=\diag(\pi,\pi^3,1)$.

We next see that $\tau$ is \emph{not} liftable.  The linear system for $\tau$ is
\begin{equation}
\begin{split}
x_1-x_2&=\alpha_{12}-\alpha_{21}=-1\\
x_1-x_3&=\alpha_{13}-\alpha_{23}=0\\
x_2-x_3&=\alpha_{23}-\alpha_{13}=0.\end{split}
\end{equation}
One checks easily that this system is inconsistent, so that $\tau$ is not liftable.

If follows that $\Aut_R(\L)=\Inn(\L)\rtimes \O_\L$, where $\O_\L$ is cyclic of order $3$, generated by conjugation by $dv(\sigma)=\Bigl(\begin{smallmatrix}0&\pi&0\\0&0&\pi^3\\1&0&0\end{smallmatrix}\Bigr)$.\qed
\end{exmp}

\begin{rem} \normalfont Suppose that $\L$ is a basic, hereditary $R$-order in $M_n(K)$. Then without loss of generality we may assume that the exponent matrix for $\L$ is given by $\alpha_{ij}=0$ if and only if $i\geq j$, so that $\L$ has the form 
\[\L=\begin{pmatrix} R&R&\dots &R&R \\ P&R&\dots &R&R\\ \vdots&\vdots &\ddots &\vdots&\vdots\\P&P&\dots&P&R\end{pmatrix}.\] 
The radical of $\L$ is obtained by replacing the $R$'s in the diagonal by $P$'s.  

If we let $P_i$ denote the $i$-th column of $\L$, then the description of $\rad(\L)$ shows that $\rad(P_i)=P_{i+1}$, where the indices are taken modulo $n$.  Consequently, the link graph of $\L$ is 

\centerline{\includegraphics*{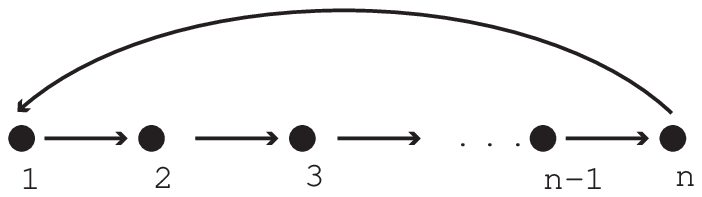}}

We see immediately that $\Aut(Q(\L))$ is cyclic of order $n$, generated by $\sigma=(12...n)$.  The corresponding linear system is $x_i-x_j=0$ for $i<j<n$, and $x_i-x_n=1$, which has the solution ${\bf x}=(0,\dots,0,-1)$.  It follows that $\sigma$ is liftable, and that $\O_\L$ is cyclic of order $n$, generated by conjugation by $\bigl(\begin{smallmatrix} 0&\pi^{-1}\\I_{n-1}&0\end{smallmatrix}\bigr)$, where $I_{n-1}$ is the $(n-1)\times (n-1)$ identity matrix.  Thus we may use Theorem \ref{main thm} to recover \cite[Proposition 4.10]{Haef Jan}\qed
\end{rem}

We call a tiled order $\L$ a \emph{$(0,1)$-order} if $\alpha_{ij}\in\{0,1\}$ for all $i,j$.  Fujita has shown that the link graph is an isomorphism invariant for basic $(0,1)$-orders \cite[Theorem]{Fujita Link}.  Using this fact, we show that every automorphism of $Q(\L)$ is liftable for such $\L$.  

\begin{thm}  Let $\L$ be a basic $(0,1)$-order.  Then every $\sigma\in\Aut(Q(\L))$ is liftable.  Consequently, $\Aut_R(\L)\cong \Inn(\L)\rtimes\Aut(Q(\L))$.
\end{thm}

\begin{pf}  Let $\sigma\in\Aut(Q(\L))$, and let $\Gamma=v(\sigma)\L v(\sigma)^{-1}$.  Since $\L$ and $\Gamma$ are isomorphic, we have $Q(\Gamma)=Q(\L)$.  By \cite[section 2]{Fujita Link}, there is a diagonal matrix $d$ such that $\Gamma=d^{-1}\Lambda d$.  It follows that 
\[\L=v(\sigma)^{-1}\Gamma v(\sigma)=v(\sigma)^{-1}d^{-1}\Lambda dv(\sigma),\]
so that $\iota_{dv(\sigma)}\in\Aut(\L)$. By Proposition \ref{liftable lemma}, $\sigma$ is liftable.\qed
\end{pf}

\section{Non-primeness of crossed products}
One of the initial motivations of this work was to extend results of \cite{Haef Jan} on the structure of crossed products over hereditary orders.  There it was shown that, if $\Gamma$ is a prime, hereditary order and $G$ is a subgroup of $\O_\Gamma$, then the crossed product $\Gamma*G$ was again prime and hereditary.  We had originally hoped that this result might generalize to the case of tiled orders as well.  However, as the following example illustrates, $\L*\O_\L$ need not be a prime order, even if $\L$ is prime and $\O_\L$ acts transitively on the primitive idempotents of $\L$.

\begin{exmp} \normalfont Let $\L=\bigl(\begin{smallmatrix} R&P\\P&R\end{smallmatrix}\bigr)$, where $R/P$ does not have characteristic $2$.  Then it is easy to see that $\O_\L$ is cyclic of order $2$, generated by conjugation by $\bigl(\begin{smallmatrix} 0&1\\1&0\end{smallmatrix}\bigr)$, which we will denote by $\tau$.  We compute directly that $\L*\O_\L\cong M_2(R)\oplus M_2(R)$ as rings, proving the claim.  

We realize $\L*\O_\L$ as a free $\L$-module of rank $2$, with basis $\{1,\bar \tau\}$, where multiplication is determined by $\bar \tau^2=1$ and $\bar \tau \lambda=\tau(\lambda)\bar\tau$ for $\lambda\in\L$.  If we let $e_{11}$ and $e_{22}$ denote the standard matrix units $\bigl(\begin{smallmatrix} 1&0\\0&0\end{smallmatrix}\bigr)$ and $\bigl(\begin{smallmatrix} 0&0\\0&1\end{smallmatrix}\bigr)$ respectively, then $\bar \tau e_{11}=e_{22}\bar\tau$ and $\bar\tau e_{22}=e_{11}\bar\tau$. If we set $e_{12}=e_{11}\bar\tau$ and $e_{21}=e_{22}\bar\tau$, then it is straightforward to check that $\{e_{11},e_{12},e_{21},e_{22}\}$ form a set of matrix units in $\L*\O_\L$. Thus $\L*\O_\L\cong M_2(S)$, where $S\cong e_1(\L*\O_\L)e_1$.  We finish by showing that $S\cong R\oplus R$ as rings.

We identify $R$ with $e_{11}\L e_{11}$ (that is, with the subset $\bigl(\begin{smallmatrix} R&0\\0&0\end{smallmatrix}\bigr)$ of $\L$), and define a homomorphism $\vphi:R[x]\rightarrow S$ by sending $1$ to $e_{11}$ and $x$ to $e_{11}\bigl(\begin{smallmatrix} 0&\pi\\\pi&0\end{smallmatrix}\bigr)\bar\tau e_{11}$.  One then computes that $x^2$ maps to $\bigl(\begin{smallmatrix} \pi^2&0\\0&0\end{smallmatrix}\bigr)\in e_{11}\L e_{11}$.  It follows that the kernel of $\vphi$ contains $x^2-\pi^2$; the fact that $S$ and $R[x]/(x^2-\pi^2)$ both have rank $2$ over $R$ shows that $\ker\vphi=(x^2-\pi^2)$.  Since the characteristic of $R/P$ is not $2$, $S\cong R[x]/(x^2-\pi^2)\cong R\oplus R$ as rings.\qed
\end{exmp}

\bibliographystyle{amsalpha}

\begin{thebibliography}{9}

\bibitem{Fujita rem} H. Fujita, A remark on tiled orders over a local Dedekind domain, \textit{Tsukuba J. Math.} 10:121--130 (1986).

\bibitem{Fujita Link} H. Fujita, Link graphs of tiled orders over a local Dedekind domain, \textit{Tsukuba J. Math.} 10:293--298 (1986).

\bibitem{Fuj Yosh} H. Fujita and H. Yoshimura, A crterion for isomorphic tiled orders over a local Dedekind domain, \textit{Tsukuba J. Math.} 16:107--111 (1992).

\bibitem{Haef Jan} J. Haefner and G. Janusz, Hereditary crossed products, \textit{Trans. Amer. Math. Soc.} 352:3381--3410 (2000).

\bibitem{Jat} V. A. Jategaonkar, Global dimension of tiled orders over a discrete valuation ring, \textit{Trans. Amer. Math. Soc.} 196:313--330 (1974).

\bibitem{Lambek} J. Lambek, \textit{Lectures on Rings and Modules (Second Edition)}, Chelsea Publishing, New York, 1976.

\bibitem{Muller FBN} B. J. M\"{u}ller, Localization in fully bounded noetherian rings, \textit{Pacific J. Math.} 67:233--245 (1976).

\bibitem{Wie Rog} A. Wiedemann and K. Roggenkamp, Path orders of global dimension two, \textit{J. Algebra} 80:113--133 (1983).



\end{thebibliography}

\end{document}